\documentclass[9pt,amsfonts]{article}
\usepackage{graphicx,amssymb,eucal,mathrsfs}
\usepackage{latexsym,amsmath,amscd,epsfig,epic,eepic,xcolor}
\usepackage[all]{xy}
\usepackage{wasysym,ntheorem}
\usepackage{epsf,graphicx,pgf,etex,amscd,pstricks}

\usepackage[hypertexnames=false,backref=page,pdftex,
 	pdfpagemode=UseNone,
 	breaklinks=true,
 	extension=pdf,
 	colorlinks=true,
 	linkcolor=blue,
 	citecolor=red,
 	urlcolor=blue,
 ]{hyperref}

\newcommand{\Z}{\mathbb{Z}}
\newcommand{\R}{\mathbb{R}}
\newcommand{\C}{\mathbb{C}}

\newcommand{\ve}{\varepsilon}

\newcommand{\SI}{\mathbb{S}}

\numberwithin{equation}{section}

\newtheorem{main}{Theorem}[section]

\newtheorem{prop}[main]{Proposition}
\newtheorem{lemma}[main]{Lemma}

\theoremstyle{definition}

\newtheorem{defi}[main]{Definition}

\title{A remark about weak fillings}
\date{March 2016}
\author{Pierre Py}

\begin{document}

\maketitle

\begin{abstract} Let $L$ be a closed manifold of dimension $n\ge 2$ which admits a totally real embedding into $\C^{n}$. Let $ST^{\ast}L$ be the space of rays of the cotangent bundle  $T^{\ast}L$ of $L$ and let $DT^{\ast}L$ be the unit disc bundle of $T^{\ast}L$ defined by any Riemannian metric on $L$. We observe that $ST^{\ast}L$ endowed with its standard contact structure admits weak symplectic fillings $W$ which are diffeomorphic to $DT^{\ast}L$ and for which any closed Lagrangian submanifold $N\subset W$ has the property that the map $H_{1}(N,\R)\to H_{1}(W,\R)$ has a nontrivial kernel. This relies on a variation on a theorem by Laudenbach and Sikorav. 
\end{abstract}

\section{Introduction}

A field of hyperplanes $\xi$ on a manifold $M$ is a contact structure if one can find locally a $1$-form $\alpha$ such that $\xi ={\rm Ker}(\alpha)$ and such that $d\alpha$ restricted to $\xi$ is everywhere nondegenerate. This implies that $M$ is odd dimensional; we write ${\rm dim} M=2n-1$. We always assume that the manifolds and contact structures we consider are compatibly oriented. This means that $\xi$ can be globally defined by a $1$-form $\alpha$ such that $\alpha \wedge (d\alpha)^{n-1}$ defines the positive orientation of $M$ and $(d\alpha)^{n-1}$ defines the positive orientation of $\xi$. The restriction of $d\alpha$ to $\xi$ is a symplectic form whose conformal class depends only on $\xi$; we denote this conformal class by ${\rm CS}_{\xi}$. 

There exist many notions of symplectic (or holomorphic) fillings of a contact manifold $(M,\xi)$. They are motivated by the fact that many examples of contact manifolds appear as boundaries of complex or symplectic manifolds with suitable convexity properties; although not all contact manifolds appear in this way. We start with an informal definition. A filling of $(M,\xi)$ is a compact symplectic manifold with boundary $(W,\omega)$ such that the oriented boundary of $W$ is identified with $M$ and such that the symplectic form $\omega$ is compatible with $\xi$ near the boundary of $W$. Depending on the additional properties one asks on $(W,\omega)$ and on the exact compatibility condition required between $\omega$ and $\xi$ one obtains many notion of fillings: weak fillings, strong fillings, exact fillings, Stein fillings... We refer the reader to~\cite{eg} for the beginning of the theory of fillings of contact manifolds and to~\cite{eliashberg,eliashberg2,eh,gay,mnw,nw,wendl2} for more recent references. Here, we will simply use the following definition, which is classical when $M$ has dimension $3$ and due to Massot, Niederkr\"uger and Wendl~\cite{mnw} in higher dimensions.  

\begin{defi}We say that $(W,\omega)$ is a weak filling of $(M,\xi)$ if the restriction $\omega_{\xi}$ of $\omega$ to the distribution $\xi$ is symplectic and if each element in the ray $\omega_{\xi}+{\rm CS}_{\xi}$ is symplectic. 
\end{defi}

When $M$ has dimension $3$, this definition reduces to the fact that $\omega_{\xi}$ is everywhere nondegenerate and defines the positive orientation of $\xi$. We refer the reader to~\cite{mnw} for some motivations for this definition and for more history of fillability problems in contact geometry. Here, the main observation that we will need is the following: 
\begin{center}
{\it If $M$ bounds a smooth strictly pseudoconvex domain $W$ in a K\"ahler manifold $(X,\omega)$, then $(W,\omega)$ is a weak filling of $(M,\xi)$, where $\xi$ is the contact structure given by the complex tangencies.}
\end{center}
We shall use this observation to produce exotic examples of weak fillings. The idea of considering open sets in K\"ahler manifolds which are strictly pseudoconvex but not {\it symplectically convex} is actually not new. It was already used by Eliashberg and Gromov 25 years ago to motivate the various definitions of convexity in complex and symplectic geometry, see~\cite[\S 3.1.3]{eg}.  
  
Consider now a manifold $L$ and its cotangent bundle $T^{\ast}L$, endowed with its Liouville form $\lambda =pdq$. Let $Y$ be the Liouville vector field of $T^{\ast}L$, whose flow is given by $(q,p)\mapsto (q,e^{t}p)$. The space $ST^{\ast}L$ of rays of $T^{\ast}L$ carries a canonical contact structure $\xi_{st}$ which can be defined as follows: identify $ST^{\ast}L$ with a hypersurface of $T^{\ast}L$ transverse to $Y$ and take the kernel of the restriction of the Liouville form to $ST^{\ast}L$. The resulting hyperplane field does not depend on this identification. Finally, we will denote by $DT^{\ast}L$ the unit disc bundle inside $T^{\ast}L$, for an auxiliary Riemannian metric on $L$. We assume that ${\rm dim} L =n \ge 2$. We will prove the following result: 
    
\begin{main}\label{wf} Let $L$ be closed. Assume that $L$ admits a totally real embedding into $\C^{n}$. Then there exists an exact symplectic form $\Omega$ on $W=DT^{\ast}L$ such that $(W,\Omega)$ is a weak filling of $(ST^{\ast}L,\xi_{st})$ and such that for any closed Lagrangian submanifold $N\subset W$, the induced map $H_{1}(N,\R)\to H_{1}(W,\R)$ has a nontrivial kernel. 
\end{main}

This applies for instance when $L$ is a torus or any oriented closed $3$-manifold. See~\cite{audin} for a discussion of which manifolds admit totally real embeddings into $\C^{n}$ and for more examples. 

When $L=T^{2}$, weak fillings of $(ST^{\ast}T^{2},\xi_{st})$ which are different from $DT^{\ast}T^{2}$ with its standard symplectic structure were already known, see for instance~\cite{giroux}. However the examples in~\cite{giroux} were not cohomologically exact, as opposed to the examples provided by Theorem~\ref{wf}. We aslo mention that Wendl~\cite{wendl} proved that every {\it exact} filling $(W,\omega)$ of $(ST^{\ast}T^{2},\xi_{st})$ is symplectomorphic to a star-shaped open set of $T^{\ast}T^{2}$. See~\cite{wendl} for the definition of exact fillings. In particular, exact fillings of $(ST^{\ast}T^{2},\xi_{st})$ are diffeomorphic to $DT^{\ast}T^{2}\simeq T^{2}\times D^{2}$ and contain an incompressible Lagrangian torus.

In cases where it is known that there is a unique exact filling of $(ST^{\ast}L,\xi_{st})$, the weak filling $(W,\Omega)$ appearing in Theorem~\ref{wf} can be symplectically embedded in $T^{\ast}L$ with its standard symplectic structure, as a deformation of the unit disc bundle. But the image of such an embedding will never contain the zero section, by our result. This embedding result follows from the fact that one can glue a product to the boundary of $W$ to obtain a manifold
$$W'=W\cup ST^{\ast}L\times [0,1]$$
such that $W'$ carries a symplectic structure which extends $\Omega$ and which makes $W'$ an exact filling of $(ST^{\ast}L,\xi_{st})$. The proof is essentially folkloric and relies on an argument of Eliashberg~\cite[Prop. 3.1]{eliashberg0}, see also~\cite[\S 2]{mnw} for more details. 

%

Let us now discuss the proof of Theorem~\ref{wf}. The weak filling that we will consider will be small tubular neighborhoods of totally real submanifolds of the complex Euclidean space. We thus consider the space $\C^{n}$ endowed with its standard symplectic form:
$$\omega_{0}=\frac{i}{2}\sum_{j=1}^{n}dz_{j}\wedge d\overline{z_{j}}.$$ 
Denote by $J_{0}$ the complex structure on $\C^{n}$. Fix a totally real embedding $$j : L \to \C^{n},$$
of a closed manifold $L$. We assume that ${\rm dim} L = n\ge 2$. If $j$ happens to be Lagrangian, we can always perturb it to a non-Lagrangian embedding, keeping it totally real. So from now on we assume that $j$ is {\it not} Lagrangian. Let $f : \C^{n}\to \R_{+}$ be the square of the distance to $j(L)$, i.e. for $p\in \C^{n}$ let
$$f(p)=\underset{q\in j(L)}{{\rm inf}}d(p,q)^{2},$$
where $d$ is the Euclidean distance in $\C^{n}$. The function $f$ is smooth and strictly plurisubharmonic near $j(L)$, see~\cite[\S 2.7]{ce}. Let
$$V_{\ve}:=\{p\in \C^{n}, f(p)\le \ve\}.$$
In the following we will say that a submanifold $N$ of an ambient manifold $W$ is $H_{1}$-embedded if the map $H_{1}(N,\R)\to H_{1}(W,\R)$ induced by the inclusion is injective. We have:

\begin{main}\label{tt} If $\varepsilon$ is small enough, the symplectic manifold $(V_{\varepsilon},\omega_{0})$ does not contain any $H_{1}$-embedded closed Lagrangian submanifold. 
\end{main}

This result is an imediate application or rather variation on a theorem by Laudenbach and Sikorav~\cite{ls}. We will prove it in section~\ref{proof}. To deduce Theorem~\ref{wf} from Theorem~\ref{tt}, we simply make the following observations. Let $\xi_{\varepsilon}$ denote the field of complex tangencies on the boundary $M_{\varepsilon}$ of $V_{\varepsilon}$:
$$\xi_{\varepsilon}(p)=T_{p}M_{\varepsilon}\cap J_{0}T_{p}M_{\varepsilon}.$$
Then, for small enough $\varepsilon$, $(M_{\varepsilon},\xi_{\varepsilon})$ is a contact manifold, $(V_{\varepsilon},\omega_{0})$ is a weak filling of $(M_{\varepsilon},\xi_{\varepsilon})$ and $M_{\varepsilon}$ and $V_{\varepsilon}$ are respectively diffeomorphic to $ST^{\ast}L$ and $DT^{\ast}L$. The only thing left to conclude the proof of Theorem~\ref{wf} is to observe that $(M_{\varepsilon},\xi_{\varepsilon})$ and $(ST^{\ast}L,\xi_{st})$ are actually contactomorphic. This is certainly well-known to experts, but we include a proof in section~\ref{ics}, due to a lack of a convenient reference. 

{\bf Acknowledgements.} I would like to thank Klaus Niederkr\"uger for his help during the preparation of this text, as well as Patrick Massot and the referee for their comments. This note was conceived in the fall of 2014, during a visit to the University of Chicago that I would like to thank for its hospitality.   

\section{Displaceability and homologically essential Lagrangians}\label{proof}

In the following, for each subset $X\subset \C^{n}$ and each positive number $\varepsilon$, we will denote by $V_{\varepsilon}(X)$ the set of points of $\C^{n}$ at distance at most $\sqrt{\varepsilon}$ from $X$. 

We start by recalling the theorem by Laudenbach and Sikorav~\cite{ls} alluded to earlier. Let $P$ be a closed manifold. Then: 
\begin{main} \cite{ls}
{\it If a sequence of smooth Lagrangian embeddings $\varphi_{\ell} : P \to \C^{n}$ $C^{0}$-converges to a smooth embedding $\varphi_{\infty} : P \to \C^{n}$, then $\varphi_{\infty}(P)$ is still Lagrangian.}
\end{main}
This can be seen as a generalization of the famous result by Eliashberg and Gromov about $C^{0}$-limits of symplectic diffeomorphisms; it also holds in more general symplectic manifolds under an additional hypothesis, see~\cite{ls}. Now let $\varphi_{\ell}$ and $\varphi_{\infty}$ be embeddings as in the statement of the Theorem. The $C^{0}$-convergence of $\varphi_{\ell}$ to $\varphi_{\infty}$ implies that for $\varepsilon$ small enough and for $\ell$ large enough, the inclusion
$$\varphi_{\ell}(P)\hookrightarrow V_{\varepsilon}(\varphi_{\infty}(P))$$
induces an isomorphism on fundamental groups. In particular it induces an injection on the first homology groups. We will see by closely inspecting Laudenbach and Sikorav's proof that this property, namely the presence of $H_{1}$-embedded Lagrangians in arbitrarily small neighborhoods of $\varphi_{\infty}(P)$ is sufficient to imply that $\varphi_{\infty}(P)$ is itself Lagrangian. It is not even necessary to assume that these Lagrangians are diffeomorphic to $P$. In other words, we will prove:

\begin{prop}\label{lsrev} Let $P\subset \C^{n}$ be a closed submanifold of real dimension $n$. Assume that there exists a sequence $(\varepsilon_{k})$ of positive numbers converging to $0$ such that for all $k$, $V_{\varepsilon_{k}}(P)$ contains a closed $H_{1}$-embedded Lagrangian submanifold. Then $P$ is Lagrangian. 
\end{prop}
 
It is clear that this proposition implies Theorem~\ref{tt} since the embedding $j$ from the previous section was assumed to be totally real but not Lagrangian. 

\medskip

We now turn to the proof of Proposition~\ref{lsrev}. We follow the proof of Laudenbach and Sikorav and check that it still applies under our hypothesis. We will use the following result of Gromov~\cite{gromov}, see also Proposition 1.2 in~\cite{sikorav}: 
\smallskip

{\it If a closed Lagrangian submanifold $N\subset \C^{n}$ is contained in $\C^{n-1}\times B(\alpha)$ then there exists a nonconstant holomorphic disc with boundary on $N$ of area smaller than $\pi \alpha^{2}$. In particular there is a loop $\gamma$ in $N$ such that $0<\langle [\lambda_{0}],\gamma\rangle <\pi \alpha^{2}$ where $[\lambda_{0}]$ is the Liouville class of $N$. Here $B(\alpha)$ stands for the ball of radius $\alpha$ in the complex plane.} 

\smallskip

We argue by contradiction and assume that $P$ is not Lagrangian. Exactly as in~\cite{ls}, we can assume that the normal bundle of $P$ in $\C^{n}$ has a nowhere vanishing section (up to replacing $P$ by $P\times  S^{1}$ in $\C^{n+1}$). We can also assume that the sequence $(\varepsilon_{k})$ is decreasing. Let $L_{k}$ be a closed $H_{1}$-embedded Lagrangian submanifold of $V_{\varepsilon_{k}}(P)$. By Theorem 1 in~\cite{ls} (see also~\cite{pol2}), there exists a Hamiltonian function $H$ with flow $\phi^{t}_{H}$ such that $\phi_{H}^{t}(P)\cap P$ is empty for all $t\in (0,\delta)$ for some positive number $\delta$. We fix a sequence $t_{\ell}$ of positive real numbers converging to $0$ and for each $\ell$ we choose $\varepsilon_{k_{\ell}}>0$ such that
\begin{equation}\label{disj}
\phi_{H}^{t_{\ell}}(V_{\ve_{k_{\ell}}}(P))\cap V_{\ve_{k_{\ell}}}(P)=\emptyset.
\end{equation}
In particular the isotopy $(\phi_{H}^{t})_{0\le t \le t_{\ell}}$ displaces the Lagrangian submanifold $L_{k_{\ell}}$. By reparametrizing it, one can find a Hamiltonian isotopy $(\varphi^{t})_{0\le t \le t_{\ell}}$ which still satisfies~\eqref{disj} and which has the additional property that it is constant equal to the identity (resp. $\varphi^{t_{\ell}}$) for $t$ close to $0$ (resp. close to $t_{\ell}$). From this fact, one can construct a particular Lagrangian embedding of $L_{k_{\ell}}\times \SI^{1}$ in $\C^{n+1}$, namely we have:

\begin{prop}\label{suspension} Identify the circle $\SI^{1}$ with $\R/2\Z$. There exists a Lagrangian embedding $\Phi_{\ell} : L_{k_{\ell}} \times \SI^{1} \to \C^{n}\times \C$ with the following properties:
\begin{enumerate}
\item the Liouville class of $\Phi_{\ell}$ is of the form $([\lambda_{\ell}],0)\in H^{1}(L_{k_{\ell}}\times \SI^{1},\mathbb{R})$,\label{propl}
\item if $\underline{\pi}$ denotes the natural projection $\C^{n}\times \C \to \C^{n}$, the map $\underline{\pi} \circ \Phi_{\ell} : L_{k_{\ell}}\times \SI^{1}\to \C^{n}$ is given by $\underline{\pi}\circ \Phi_{\ell}(x,t)=\varphi^{t\cdot t_{\ell}}(x)$ for $t\in [0,1]$ and $\underline{\pi} \circ \Phi_{\ell}(x,t)=\varphi^{(2-t)\cdot t_{\ell}}(x)$ for $t\in [1,2]$,\label{ppp} 
\item the image of $\Phi_{\ell}$ is contained in $\C^{n}\times B(\alpha_{\ell})$ where $\alpha_{\ell}$ goes to $0$ as $\ell$ goes to infinity. 
\end{enumerate} 
\end{prop}
This proposition is very classical, we refer the reader to~\cite[2.3.$B_{3}'$]{gromov} or~\cite{pol} for its proof. Exactly as in~\cite{ls}, we will now use Proposition~\ref{suspension} to finish the proof of Proposition~\ref{lsrev}. According to Gromov's result mentioned above, there exists a non-constant holomorphic disc $g^{\ell} =(g^{\ell}_{1},g^{\ell}_{2}) : D \to \C^{n}\times \C$ with boundary contained in $\Phi_{\ell}(L_{k_{\ell}} \times \SI^{1})$ of area $\le \pi \alpha_{\ell}^{2}$. Let $p : L_{k_{\ell}}\times \SI^{1} \to L_{k_{\ell}}$ be the natural projection. The map $$h_{\ell}:=p\circ \Phi_{\ell}^{-1}\circ g^{\ell} : \partial D \to L_{k_{\ell}}$$ represents a nontrivial homology class by Proposition~\ref{suspension}~\eqref{propl}.  Recall that the sequence $(\ve_{k})$ is decreasing, hence we can consider the map 
$$i \circ h_{\ell} : \partial D \to V_{\ve_{1}}(P),$$
where $i : L_{k_{\ell}} \to V_{\ve_{1}}(P)$ is the inclusion. It still represents a nontrivial homology class since $L_{k_{\ell}}$ is $H_{1}$-embedded. 

\begin{lemma} The maps $i\circ h_{\ell}$ and $g^{\ell}_{1}$, considered as maps from $\partial D$ to $V_{\ve_{1}}(P)$ are homotopic.
\end{lemma}
\noindent {\it Proof.} Define $v_{\ell} : \partial D \to L_{k_{\ell}}\times \SI^{1}$ by $v_{\ell}=\Phi_{\ell}^{-1}\circ g^{\ell}$. We have $i \circ h_{\ell}=i \circ p\circ v_{\ell}$ and $g^{\ell}_{1}=\underline{\pi}\circ \Phi_{\ell}\circ v_{\ell}$. Property~\eqref{ppp} from Proposition~\ref{suspension} implies that the two maps
$$i \circ p , \underline{\pi}\circ \Phi_{\ell} : L_{k_{\ell}}\times \SI^{1} \to V_{\ve_{1}}(P)$$
are homotopic, hence $i\circ h_{\ell}$ and $g^{\ell}_{1}$ are homotopic.\hfill $\Box$

Hence $g^{\ell}_{1}(\partial D)$ represents a nontrivial homology class in $H_{1}(V_{\ve_{1}}(P),\R)$. However, $g^{\ell}_{1}(\partial D)$ bounds the disc $g^{\ell}_{1}(D)$ whose area is bounded above by the area of $g^{\ell}(D)$, which goes to $0$ as $\ell$ goes to infinity. This gives a contradiction with the following Lemma from~\cite[p. 165]{ls}, applied to $A=V_{\ve_{1}}(P)$.  

\begin{lemma} Let $A$ be a compact domain in $\C^{n}$. Then there exists $\delta (A) >0$ such that for any smooth map $u : D \to \C^{n}$ with boundary contained in $A$, and such that $[u(\partial D)]\neq 0$ in $H_{1}(A,\R)$, one has:
$$area (D)\ge \delta (A).$$
\end{lemma}

This concludes the proof of Proposition~\ref{lsrev} and thus of Theorem~\ref{tt}. As the reader will have noticed, we have only repeated the proof of Laudenbach and Sikorav~\cite{ls}.

\section{Identification of the contact structure}\label{ics}

Recall from the introduction that $M_{\ve}$ is the boundary of the $\ve$-tubular neighborhood of a totally real submanifold $j(L)$ of $\C^{n}$, endowed with its canonical contact structure $\xi_{\ve}$. We prove here:

\begin{prop}\label{iii} The manifold $(M_{\varepsilon},\xi_{\varepsilon})$ is contactomorphic to $(ST^{\ast}L,\xi_{st})$.
\end{prop}

In what follows, we will identify $L$ and $j(L)$ and will not mention the map $j$ anymore. Hence, we will think of $TL$ as a subbundle of $L\times \C^{n}$. 

First, we will have to deal with the fact that the normal bundle $NL$ of $L$ in $\C^{n}$ need not  coincide with the image under $J_{0}$ of its tangent bundle; indeed $L$ is not assumed to be Lagrangian. For this, we choose once and for all a smooth map
$$\phi : TL\times [0,1]\to L\times \C^{n}$$
with the following properties:
\begin{itemize}
\item the map $\phi$ is of the form $\phi(x,v,t)=(x,A(x,v,t))$,
\item for each $t$, $\phi_{t}:=\phi (\cdot , \cdot , t)$ is an isomorphism onto a rank $n$ subbundle of $L\times \C^{n}$ which is transverse to $TL$,
\item $\phi_{0}(TL)=NL$, 
\item $\phi_{1}=J_{0}$ and hence $\phi_{1}(TL)=J_{0}TL$,
\item the metric on $TL$ induced by the embedding $\phi_{t}$ and by the Euclidean metric on $\C^{n}$ does not depend on $t$, it is denoted by $\vert \cdot \vert$.  
\end{itemize}
It is easy to construct a map $\phi$ satisfying the first four properties, and one can then always achieve the last one by composing with a bundle automorphism $TL\to TL$ depending on $t$. In the following we will denote by $D_{\delta}(TL)$ the open disc bundle of radius $\delta$ in $TL$, for the metric $\vert \cdot \vert $ induced by $\phi$. Consider now the maps 
$$\theta_{t} : TL\to \C^{n}\;\;\;\;\; (t\in [0,1])$$
defined by $\theta_{t}(x,v)=j(x)+A(x,v,t)$. There exists $\delta_{1}>0$ such that for each $t$ in $[0,1]$, the map $\theta_{t}$ is injective and a local diffeomorphism on $D_{\delta_{1}}(TL)$. The open set 
$$\theta_{t}(D_{\delta_{1}}(TL))\subset \C^{n}.$$
is a tubular neighborhood of $L$ modeled on a varying subbundle of $L\times \C^{n}$, which is always transverse to $TL$. Let $J^{t}$ be the almost complex structure on $D_{\delta_{1}}(TL)$ which is the pull back of $J_{0}$ by $\theta_{t}$. 
 Now recall the following formula~\cite[\S 2.2]{ce}. If $\varphi$ is a function defined on an open set $V$ of $\C^{n}$, if $x$ is a point of $V$ and $u\in \C^{n}$, we have:
 $$-dd^{\C}\varphi_{x}(u,J_{0}u)={\rm Hess}_{\varphi, x}(u)+{\rm Hess}_{\varphi , x}(J_{0}u)$$
 where ${\rm Hess}_{\varphi , x}$ is the Hessian of $\varphi$ at $x$ and $d^{\C}\varphi=d\varphi \circ J_{0}$. Let $h : D_{\delta_{1}}(TL)\to \R_{+}$ be the function $h(x,v)=\vert v\vert^{2}$. By applying the above formula to the functions $h\circ \theta_{t}^{-1}$ for $x\in L$ one finds: 
 $$-dd^{\C}_{J^{t}}h>0 \;\;\;\;\; (t\in [0,1])$$
along the zero section $L\subset D_{\delta_{1}}(TL)$. Here we have used the notation $$d^{\C}_{J^{t}}h=dh\circ J^{t}.$$
This implies that there exists a positive number $\delta_{2}<\delta_{1}$ such that 
$$-dd^{\C}_{J^{t}}h>0$$
on all of $D_{\delta_{2}}(TL)$ for all $t$ in $[0,1]$. Fix any number $\ve \in (0,\delta_{2}^{2})$ and let $$S_{\varepsilon}=\{h=\ve \}\subset D_{\delta_{2}}(TL).$$
This is the sphere bundle of radius $\sqrt{\ve}$. Let $\xi (t)=TS_{\ve}\cap J^{t}TS_{\ve}$ be the field of complex tangencies for $J^{t}$ in $S_{\ve}$. All the $\xi (t)$'s are contact structures since $h$ is $J^{t}$-convex on $D_{\delta_{2}}(TL)$. Note that $\xi (0)=\xi_{\ve}$ is the contact structure that we want to identify with the canonical contact structure on $ST^{\ast}L$. By Gray's theorem, all the $\xi (t)$'s are isomorphic, hence it is enough to prove that $(S_{\ve},\xi (1))$ is contactomorphic to $(ST^{\ast}L,\xi_{st})$.  

Recall now that the Riemannian metric $\vert \cdot \vert$ on $L$ induces a decomposition of the tangent bundle of $TL$ into horizontal and vertical subbundles. At each point $(x,v)\in TL$ one has a decomposition
$$T_{(x,v)}TL=H(x,v)\oplus V(x,v)$$
where $V(x,v)$ is the tangent space to the fiber of the projection $\pi : TTL \to L$ at $(x,v)$ and $H(x,v)$ is the horizontal subspace defined by the Levi-Civita connection of $\vert \cdot \vert$, see~\cite[\S 1.3]{paternain}. Both $V(x,v)$ and $H(x,v)$ are canonically identified with $T_{x}L$ hence one has an identification
$$T_{(x,v)}TL\simeq T_{x}L\times T_{x}L$$
(where the first factor is horizontal and the second vertical). Let $J^{\ast}$ be the almost complex structure on $TL$ defined by $(u,v)\mapsto (-v,u)$ under the previous identification (see~\cite[\S 1.3.2]{paternain} for more details). Let $\xi^{\ast}$ be the field of complex tangencies for $J^{\ast}$ on $S_{\ve}$. Then $(S_{\ve},\xi^{\ast})$ is contactomorphic to $(ST^{\ast}L,\xi_{st})$, see~\cite{paternain}. 

We now want to find a path of almost complex structures from $J^{\ast}$ to $J^{1}$, to relate the contact structures $\xi^{\ast}$ and $\xi(1)$. But along the zero section both $J^{1}$ and $J^{\ast}$ interchange the vertical and horizontal subbundles and actually the hypothesis $\phi_{1}=J_{0}$ made earlier implies that $J^{\ast}=J^{1}$ along the zero section. Since $J^{\ast}$ is tamed by a symplectic form on $TL$, this implies that $J^{1}$ and $J^{\ast}$ are tamed by a common symplectic form, say $\Omega$, on $D_{\delta_{3}}(TL)$ for $\delta_{3}\in (0,\delta_{2})$ small enough. Now we can find a path $J(t)$ of almost complex structures on $D_{\delta_{3}}(TL)$ such that $J(0)=J^{\ast}$ and $J(1)=J^{1}$ and such that $J(t)$ does not depend on $t$ along the zero section. Since 
$$-dd^{\C}_{J(t)}h>0$$
along the zero section for all $t$, we can shrink once more our neighborhood and conclude that 
$$-dd^{\C}_{J(t)}h>0$$
for all $t$ on $D_{\delta_{4}}(TL)$ for some $\delta_{4}\in (0,\delta_{3})$. Appealing to Gray's theorem again and shrinking $\ve$ if necessary, we conclude that $(S_{\ve},\xi (1))$ and $(S_{\ve},\xi^{\ast})$ are contactomorphic. This concludes the proof of Proposition~\ref{iii}.


\bigskip
\bigskip

\begin{small}
\begin{tabular}{l}
IRMA, Universit\'e de Strasbourg \& CNRS\\
67084 Strasbourg, France\\
ppy@math.unistra.fr\\    
\end{tabular}
\end{small}


\begin{thebibliography}{00}

\bibitem{audin} M.~Audin, \emph{Fibr\'es normaux d'immersions en dimension double, points doubles d'immersions lagrangiennes et plongements totalement r\'eels}, Comment. Math. Helv.~{\bf 63}, No.~4 (1988), 593--623.  

\bibitem{ce} K.~Cieliebak and Y.~Eliashberg, \emph{From Stein to Weinstein and back. Symplectic geometry of affine complex manifolds.}, American Mathematical Society Colloquium Publications~{\bf 59}, American Mathematical Society, Providence, RI (2012).    

\bibitem{eliashberg0} Y.~Eliashberg, \emph{On symplectic manifolds with some contact properties}, J. Differential Geom.~{\bf 33}, No.~1 (1991), 233--238. 

\bibitem{eliashberg} Y.~Eliashberg, \emph{Unique holomorphically fillable contact structure on the $3$-torus}, Internat. Math. Res. Notices No.~2 (1996), 77--82. 

\bibitem{eliashberg2} Y.~Eliashberg, \emph{A few remarks about symplectic fillings}, Geom. Topol.~{\bf 8} (2004), 277--293. 

\bibitem{eg} Y.~Eliashberg and M.~Gromov, \emph{Convex symplectic manifolds}, Several complex variables and complex geometry, Part 2 (Santa Cruz, CA, 1989), Proc. Sympos. Pure Math.~{\bf 52}, Amer. Math. Soc., Providence, RI (1991), 135--162.  

\bibitem{eh} J. B.~Etnyre and K.~Honda, \emph{Tight contact structures with no symplectic fillings}, Invent. Math.~{\bf 148}, No.~3 (2002), 609--626. 

\bibitem{gay} D. T.~Gay, \emph{Four-dimensional symplectic cobordisms containing three-handles}, Geom. Topol.~{\bf 10} (2006), 1749--1759. 

\bibitem{giroux} E.~Giroux, \emph{Une structure de contact, m\^eme tendue, est plus ou moins tordue}, Ann. Sci. \'Ecole Norm. Sup. (4)~{\bf 27}, No.~6 (1994), 697--705. 

\bibitem{gromov} M.~Gromov, \emph{Pseudoholomorphic curves in symplectic manifolds}, Invent. Math.~{\bf 82}, No.~2 (1985), 307--347.  

\bibitem{ls} F.~Laudenbach and J.-C.~Sikorav, \emph{Hamiltonian disjunction and limits of Lagrangian submanifolds}, Internat. Math. Res. Notices No.~4 (1994).

\bibitem{mnw} P.~Massot, K.~Niederkr\"uger and C.~Wendl, \emph{Weak and strong fillability of higher dimensional contact manifolds}, Invent. Math.~{\bf 192}, No.~2 (2013), 287--373.  

\bibitem{nw} K.~Niederkr\"uger and C.~Wendl, \emph{Weak symplectic fillings and holomorphic curves}, Ann. Sci. \'Ec. Norm. Sup\'er. (4) {\bf 44}, No.~5 (2011), 801--853.

\bibitem{paternain} G. P.~Paternain, \emph{Geodesic Flows}, Progress in Mathematics~{\bf 180}, Birkha\"user Boston, Inc., Boston, MA (1999). 

\bibitem{pol} L.~Polterovich, \emph{Symplectic displacement energy for Lagrangian submanifolds}, Ergodic Theory and Dynamical Systems~{\bf 13}, No.~2 (1993), 357--367. 

\bibitem{pol2} L.~Polterovich, \emph{An obstacle to non-Lagrangian intersections}, The Floer memorial volume, Progr. Math.~{\bf 133}, Birkha\"user, Basel (1995), 575--586. 

\bibitem{sikorav} J.-C.~Sikorav, \emph{Quelques propri\'et\'es des plongements lagrangiens}, M\'em. Soc. Math. France (N. S.)~{\bf 46} (1991), 151--167. 

\bibitem{wendl} C.~Wendl, \emph{Strongly fillable contact manifolds and $J$-holomorphic foliations}, Duke Math. J.~{\bf 151}, No.~3 (2010), 337--384. 

\bibitem{wendl2} C.~Wendl, \emph{A hierarchy of local symplectic filling obstructions for contact 3-manifolds}, Duke Math. J.~{\bf 162}, No.~12 (2013), 2197--2283. 

\end{thebibliography}
\end{document}